\documentclass[14pt]{article}
\usepackage{mathrsfs}
\usepackage{amsthm}
\usepackage{amsmath}
\usepackage{graphicx}
\usepackage{color}
\usepackage{amsfonts}
\newtheorem{theorem}{Theorem}[section]

\newtheorem{lemma}[theorem]{Lemma}

\numberwithin{equation}{section}
\normalsize

\begin{document}
\title{\textbf{Mean field limit for survival probability of the high-dimensional contact process}}

\author{Xiaofeng Xue \thanks{\textbf{E-mail}: xfxue@bjtu.edu.cn \textbf{Address}: School of Science, Beijing Jiaotong University, Beijing 100044, China.}\\ Beijing Jiaotong University}

\date{}
\maketitle

\noindent {\bf Abstract:}
In this paper we are concerned with the contact process on the squared lattice. The contact process intuitively describes the spread of the infectious disease on a graph, where an infectious vertex becomes healthy at a constant rate while a healthy vertex is infected at rate proportional to the number of infectious neighbors. As the dimension of the lattice grows to infinity, we give a mean field limit for the survival probability of the process conditioned on the event that only the origin of the lattice is infected at $t=0$. The binary contact path process is a main auxiliary tool for our proof.


\section{Introduction and main result}\label{section 1}
In this paper we are concerned with the contact process on the squared lattice $\mathbb{Z}^d$. For any $x,y\in \mathbb{Z}^d$, we write $x\sim y$ when and only when there is an edge connecting $x$ and $y$, i. e. $x$ and $y$ are neighbors. For later use, we denote by $O$ the origin of $\mathbb{Z}^d$.

The contact process on $\mathbb{Z}^d$ is a continuous-time Markov process with state space $\{0,1\}^{\mathbb{Z}^d}$, that is to say, at each vertex there is a spin taking value $0$ or $1$. For any configuration $\eta\in \{0,1\}^{\mathbb{Z}^d}$ and $x\in \mathbb{Z}^d$, we denote by $\eta(x)$ the value of the spin at $x$. For any $t\geq 0$, we denote by $\eta_t$ the configuration of the contact process at the moment $t$, then $\{\eta_t\}_{t\geq 0}$ evolves according to the following transition rates function. For each $x\in \mathbb{Z}^d$,
\begin{equation*}
\eta(x)\rightarrow 1-\eta(x)\text{~at rate~}
\begin{cases}
1&\text{~if\quad} \eta(x)=1,\\
\frac{\lambda}{2d}\sum\limits_{y:y\sim x}\eta(y) &\text{~if\quad} \eta(x)=0,
\end{cases}
\end{equation*}
where $\lambda$ is a positive constant called the infection rate and $d$ is the dimension of the lattice.

The contact process intuitively describes the spread of an infectious disease on a graph. Vertices in state $1$ are infectious while vertices in state $0$ are healthy. An infectious vertex becomes healthy at rate one while an healthy vertex is infected at rate proportional to the number of infectious neighbors.

The contact process is first introduced in \cite {Har1974} by Harris in 1974 and has been an important model for the development of the theory of interacting particle systems since then. For a detailed survey of the study of the contact process, see Chapter Six of \cite{Lig1985} and Part one of \cite{Lig1999}.

To give our main result, we introduce some notations and definitions. For $\lambda>0$ and $d\geq 1$, we denote by $P_{\lambda,d}$ the probability measure of the contact process $\{\eta_t\}_{t\geq 0}$ on $\mathbb{Z}^d$ with infection rate $\lambda$. For later use, for any $t>0$, we identify $\eta_t$ with the set
\[
\{x\in \mathbb{Z}^d:~\eta_t(x)=1\},
\]
which is the set of vertices in state $1$ at the moment $t$. For any $A\subseteq \mathbb{Z}^d$, we write $\eta_t$ as $\eta_t^A$ when $\eta_0=A$. When $A=\{x\}$ for some $x\in \mathbb{Z}^d$, we write $\eta_t^A$ as $\eta_t^x$ instead of $\eta_t^{\{x\}}$.

Our main result gives a mean field limit for the survival probability of the contact process conditioned on $\eta_0=\{O\}$ as the dimension $d$ grows to infinity.
\begin{theorem}\label{theorem main 1.1}
\begin{equation}\label{equ 1.1}
\lim_{d\rightarrow+\infty}P_{\lambda,d}(\eta_t^O\neq \emptyset,\forall~t\geq 0)=
\begin{cases}
\frac{\lambda-1}{\lambda} &\text{~if\quad} \lambda\geq 1,\\
0 & \text{~if\quad}\lambda<1.
\end{cases}
\end{equation}
\end{theorem}

\proof[Remark 1] It is obviously that $P_{\lambda,d}(\eta_t^{O}\neq \emptyset,\forall~t\geq 0)$ is increasing with the infection rate $\lambda$, hence it is reasonable to define the following critical value
\[
\lambda_c(d)=\sup\{\lambda:~P_{\lambda,d}(\eta_t^O\neq \emptyset,\forall~t\geq 0)=0\}
\]
for each $d\geq 1$. $\lambda_c$ is called the critical value of the contact process. It is shown in \cite{Grif1983} and \cite{Hol1981} that
\begin{equation}\label{equ 1.2}
\lim_{d\rightarrow+\infty}\lambda_c(d)=1.
\end{equation}
The case $\lambda<1$ for our main result follows from Equation \eqref{equ 1.2} directly, since for any $\lambda<1$, $\lambda<\lambda_c(d)$ for sufficiently large $d$. For the case $\lambda>1$, it is shown in \cite{Grif1983} that
\[
\liminf_{d\rightarrow+\infty}P_{\lambda,d}(\eta_t^O\neq \emptyset,\forall~t\geq 0)\geq \frac{\lambda-1}{2\lambda}.
\]
Our result shows that the factor $1/2$ is not necessary.

\qed

\proof[Remark 2]

Our result is called the mean field limit because we deduce it according to a mean field ODE. By the duality relationship between the contact process and itself
(see Section 3.4 of \cite{Lig1985}),
\[
P_{\lambda,d}(\eta_t^O\neq \emptyset)=P_{\lambda,d}(\eta_t^{\mathbb{Z}^d}(O)=1)
\]
and hence
\[
P_{\lambda,d}(\eta_t^O\neq \emptyset,\forall~t\geq 0)=\lim_{t\rightarrow +\infty}P_{\lambda,d}(\eta_t^{\mathbb{Z}^d}(O)=1).
\]
According to Hille-Yosida Theorem and the transition rates function of the contact process,
\begin{align*}
\frac{d}{dt}P_{\lambda,d}(\eta_t^{\mathbb{Z}^d}(O)=1)=&-P_{\lambda,d}(\eta_t^{\mathbb{Z}^d}(O)=1)\\
&+\frac{\lambda}{2d}\sum_{x:x\sim O}P_{\lambda,d}(\eta_t^{\mathbb{Z}^d}(O)=0,~\eta_t^{\mathbb{Z}^d}(x)=1).
\end{align*}
For large $d$, according to the mean field analysis, we assume that $\eta_t(O)$ and $\eta_t(x)$ where $x\sim O$ are approximately independent with each other, then $P_{\lambda,d}(\eta_t^{\mathbb{Z}^d}(O)=1)$ approximately equals to the solution of the ODE
\[
\begin{cases}
&\frac{d}{dt}f_t=-f_t+\lambda f_t(1-f_t),\\
&f_0=1.
\end{cases}
\]
When $\lambda>1$, the above ODE has a positive fix point $f_\infty=\frac{\lambda-1}{\lambda}$ while when $\lambda\leq 1$, the solution of the ODE converges to $0$ as $t\rightarrow+\infty$.

\qed

The proof of Theorem \ref{theorem main 1.1} is given in the next section. The binary contact path process introduced in \cite{Grif1983} is a main auxiliary tool for our proof.

\section{Proof of Theorem \ref{theorem main 1.1}}\label{section 2}
In this section we give the proof of Theorem \ref{theorem main 1.1}. For later use, we introduce some notations and definitions first. For each $d\geq 1$, we denote by $\{S_n^{(d)}:n=0,1,2,\ldots\}$ the discrete-time simple random walk on $\mathbb{Z}^d$. We define
\[
\sigma(d)=\inf\{n\geq 0:~S_n^{(d)}=O\}
\]
as the first moment that $\{S_n^{(d)}\}_{n\geq 0}$ visits the origin $O$. For $1\leq i \leq d$, we use $e_i(d)$ to denote the $i$th elementary unit vector of $\mathbb{Z}^d$, i. e.
\[
e_i(d)=(0,\ldots,0,\underbrace{1}_{i \text{th}},0,\ldots,0).
\]
For each $d\geq 1$, we define
\[
H(d)=P\Big(\sigma (d)<+\infty\Big|S_0^{(d)}=e_1(d)\Big)
\]
as the probability that $O$ is visited at least once conditioned on $S_0^{(d)}=e_1(d)$, then the following lemma is crucial for us to prove Theorem \ref{theorem main 1.1}.

\begin{lemma}\label{lemma 2.1}
For $\lambda>1$ and any finite nonempty set $A\subseteq \mathbb{Z}^d$,
\begin{align}\label{equ 2.1}
&P_{\lambda,d}(\eta_t^A\neq \emptyset,\forall~t\geq 0)\notag\\
&\geq \frac{|A|^2\big(\lambda-1-2\lambda H(d)\big)}{(|A|^2-|A|)(\lambda-1)\big(1-H(d)\big)+2|A|\lambda\big(1-H(d)\big)},
\end{align}
where $|A|$ is the cardinality of the set $A$.
\end{lemma}

We give the proof of Lemma \ref{lemma 2.1} at the end of this section. Now we show how to utilize Lemma \ref{lemma 2.1} to prove Theorem \ref{theorem main 1.1}.

\proof[Proof of Theorem \ref{theorem main 1.1}]

We only need to deal with the case where $\lambda>1$ according to Remark 1 given in Section \ref{section 1}. For the contact process $\{\eta_t^O\}_{t\geq 0}$ on $\mathbb{Z}^d$, we use $\{Y_n\}_{n\geq 0}$ to denote the embedded chain of $\{|\eta_t^O|\}_{t\geq 0}$, where $|\eta_t^O|$ is the number of vertices in state $1$ at the moment $t$. According to the transition rates function of the contact process, $|\eta_t|$ flips to $|\eta_t|-1$ at rate $|\eta_t|$ while $|\eta_t|$ flips to $|\eta_t|+1$ at rate at most
\[
2d|\eta_t|\frac{\lambda}{2d}=\lambda|\eta_t|,
\]
since each infectious vertex has at most $2d$ healthy neighbors to infect. As a result, $\{Y_n\}_{n\geq 0}$ is stochastically dominated from above by the asymmetric random walk $\{X_n\}_{n\geq 0}$ on $\mathbb{Z}^1$ such that
\[
P(X_{n+1}-X_n=1)=\frac{\lambda}{\lambda+1}=1-P(X_{n+1}-X_n=-1).
\]
As a result, according to classic theory of the random walk on $\mathbb{Z}_1$,
\begin{align*}
P_{\lambda,d}(\eta_t^O\neq \emptyset,\forall~t\geq 0)&=P\big(Y_n\neq 0,\forall~n\geq 1\big|Y_0=1\big)\\
&\leq P\big(X_n\neq 0,\forall~n\geq 1\big|X_0=1\big)=\frac{\lambda-1}{\lambda}
\end{align*}
for $\lambda>1$ and hence
\begin{equation}\label{equ 2.2}
\limsup_{d\rightarrow+\infty}P_{\lambda,d}(\eta_t^O\neq \emptyset,\forall~t\geq 0)\leq \frac{\lambda-1}{\lambda}.
\end{equation}

For any integer $K$, we define $\theta_K=\inf\{t:~|\eta_t^O|=K\}$ as the first moment that the number of infectious vertices is $K$ and define
$\tau_K=\inf\{n\geq 1:~Y_n=K\}$, then for each $K\geq 1$,
\begin{equation*}
P_{\lambda,d}(\theta_K<+\infty)=P\big(\tau_K<\tau_0\big|Y_0=1\big).
\end{equation*}
For $t<\theta_K$, $|\eta_t|$ flips to $|\eta_t|-1$ at rate $|\eta_t|$ while $|\eta_t|$ flips to $|\eta_t|+1$ at rate at least $\frac{\lambda}{2d}(2d-K)|\eta_t|=\lambda(1-\frac{K}{2d})|\eta_t|$, since each infectious vertex has at least $(2d-K)$ healthy neighbors. As a result, for $n<\tau_0\wedge \tau_K$, $\{Y_n\}_{n<\tau_0\wedge\tau_K}$ is stochastically dominated from below by the asymmetric random walk $\{V_n\}_{n\geq 0}$ on $\mathbb{Z}^1$ such that
\[
P(V_{n+1}-V_n=1)=\frac{\lambda(1-\frac{K}{2d})}{1+\lambda(1-\frac{K}{2d})}=1-P(V_{n+1}-V_n=-1).
\]
Let $\beta_K=\inf\{n\geq 1:~V_n=K\}$, then according to the classic theory of the random walk on $\mathbb{Z}^1$,
\begin{align*}
P_{\lambda,d}(\theta_K<+\infty)&=P\big(\tau_K<\tau_0\big|Y_0=1\big)\\
&\geq P\big(\beta_K<\beta_0\big|V_0=1\big)=\frac{1-\frac{1}{\lambda(1-\frac{K}{2d})}}{1-\big(\frac{1}{\lambda(1-\frac{K}{2d})}\big)^K}.
\end{align*}
As a result, by Lemma \ref{lemma 2.1} and the strong Markov property,
\begin{align*}
&P_{\lambda,d}(\eta_t^O\neq \emptyset,\forall~t\geq 0)\\
&\geq P_{\lambda,d}(\theta_K<+\infty)\frac{K^2\big(\lambda-1-2\lambda H(d)\big)}{(K^2-K)(\lambda-1)\big(1-H(d)\big)+2K\lambda\big(1-H(d)\big)}\\
&\geq \frac{1-\frac{1}{\lambda(1-\frac{K}{2d})}}{1-\big(\frac{1}{\lambda(1-\frac{K}{2d})}\big)^K}\frac{K^2\big(\lambda-1-2\lambda H(d)\big)}{(K^2-K)(\lambda-1)\big(1-H(d)\big)+2K\lambda\big(1-H(d)\big)}.
\end{align*}
It is shown in \cite{Kesten1964} that $\lim_{d\rightarrow+\infty}H(d)=0$ (\cite{Kesten1964} gives a more precise result that $H(d)\sim 1/(2d)$), hence
\[
\liminf_{d\rightarrow+\infty}P_{\lambda,d}(\eta_t^O\neq \emptyset,\forall~t\geq 0)\geq \frac{1-\frac{1}{\lambda}}{1-(\frac{1}{\lambda})^K}\frac{K^2(\lambda-1)}{(K^2-K)(\lambda-1)+2K\lambda}
\]
for each $K\geq 1$. Let $K\rightarrow+\infty$, then we have
\begin{equation}\label{equ 2.3}
\liminf_{d\rightarrow+\infty}P_{\lambda,d}(\eta_t^O\neq \emptyset,\forall~t\geq 0)\geq \frac{\lambda-1}{\lambda}.
\end{equation}
Theorem \ref{theorem main 1.1} follows from Equations \eqref{equ 2.2} and \eqref{equ 2.3} directly.

\qed

At last we only need to give the proof of Lemma \ref{lemma 2.1}. The binary contact path process $\{\zeta_t\}_{t\geq 0}$ introduced in \cite{Grif1983} is utilized in our proof. The state space of $\zeta_t$ is $[0,+\infty)^{\mathbb{Z}^d}$, i. e. at each vertex there is a spin taking a non-negative value. For each $x\in \mathbb{Z}^d$, $\zeta(x)$ flips to $0$ at rate $1$ while flips to $\zeta(x)+\zeta(y)$ at rate $\frac{\lambda}{2d}$ for each neighbor $y$ of $x$. Between the moments that $\zeta(x)$ flips, $\zeta(x)$ evolves according to the ODE
\[
\frac{d}{dt}\zeta_t(x)=(1-\lambda)\zeta_t(x).
\]
We assume that $\zeta_0(x)=1$ for each $x\in \mathbb{Z}^d$, then it is easy to see that
\begin{equation}\label{equ 2.3 two}
\eta_t^{\mathbb{Z}^d}=\{x:~\zeta_t(x)>0\}
\end{equation}
in the sense of coupling. An intuitive explanation of Equation \eqref{equ 2.3 two} is that for the binary contact path process we consider the seriousness of the disease of the infectious vertex. An infectious vertex can be further infected by being added the seriousness of the disease of neighbors. If we distinguish the vertices according to whether they are infectious, then we obtain the contact process. Now we can give the proof of Lemma \ref{lemma 2.1}.

\proof[Proof of Lemma \ref{lemma 2.1}]

According to the duality relationship between the contact process and itself (See Section 3.4 of \cite{Lig1985}),
\begin{equation*}
P_{\lambda,d}(\eta_t^A\neq \emptyset)=P_{\lambda,d}(\eta_t^{\mathbb{Z}^d}(x)=1 \text{~for some~}x\in A).
\end{equation*}
and hence
\begin{equation}\label{equ 2.4}
P_{\lambda,d}(\eta_t^A\neq \emptyset,\forall~t\geq 0)=\lim_{t\rightarrow+\infty}P_{\lambda,d}(\eta_t^{\mathbb{Z}^d}(x)=1 \text{~for some~}x\in A).
\end{equation}
By Equation \eqref{equ 2.3 two} and H\"{o}lder's inequality, utilizing the spatial homogeneity of $\{\zeta_t\}_{t\geq 0}$,
\begin{align}\label{equ 2.5}
&P_{\lambda,d}(\eta_t^{\mathbb{Z}^d}(x)=1 \text{~for some~}x\in A)=P_{\lambda,d}(\sum_{x\in A}\zeta_t(x)>0) \notag\\
&\geq \frac{\big(E\sum_{x\in A}\zeta_t(x)\big)^2}{E\Big(\big(\sum_{x\in A}\zeta_t(x)\big)^2\Big)}=\frac{|A|^2(E\zeta_t(O))^2}{\sum\limits_{x\in A}\sum\limits_{y\in A}F_t(x-y)},
\end{align}
where $E$ is the expectation operator with respect to $P_{\lambda,d}$ and
\[
F_t(u)=E(\zeta_t(O)\zeta(u))
\]
for any $u\in \mathbb{Z}^d$. In Chapter 9 of \cite{Lig1985}, Liggett extends the Hille-Yosida Theorem for the linear systems, which the binary contact path process belongs to. By utilizing the extensive version of Hille-Yosida Theorem, it is easy to see that
\begin{equation}\label{equ 2.6}
\frac{d}{dt} E\zeta_t(O)=-E\zeta_t(O)+2d\frac{\lambda}{2d}E\zeta_t(O)+(1-\lambda)E\zeta_t(O)=0
\end{equation}
while
\begin{equation}\label{equ 2.7}
\frac{d}{dt} F_t=GF_t,
\end{equation}
where $G$ is a $\mathbb{Z}^d\times \mathbb{Z}^d$ matrix such that
\[
G(x,y)=
\begin{cases}
-2\lambda &\text{~if~} x\neq O\text{~and~} y=x,\\
\frac{\lambda}{d} &\text{~if~} x\neq O\text{~and~} y\sim x,\\
1-\lambda & \text{~if~} x=y=O,\\
2\lambda &\text{~if~}x=O\text{~and~}y=e_1(d).
\end{cases}
\]
Note that the spatial homogeneity of $\zeta_t$ is utilized to obtain the above two equations. By Equation \eqref{equ 2.6},
\begin{equation}\label{equ 2.9}
E\zeta_t(O)=1
\end{equation}
for any $t\geq0$.
According to the classic theory of ODE in the Banach space, it is easy to check that ODE \eqref{equ 2.7} satisfies Lipschitz condition under a $l_\infty$ norm and hence ODE \eqref{equ 2.7} has the unique solution
\[
F_t=A_tF_0,
\]
where $A_t=e^{tG}=\sum_{n=0}^{+\infty}\frac{t^nG^n}{n!}$. It is easy to check that the sum in the definition of $A_t$ converges and $A_t(x,y)\geq 0$ for any $x,y\in \mathbb{Z}^d$ according to the definition of $G$. For any $x\in \mathbb{Z}^d$, we define
\[
L(x)=P\Big(\sigma(d)<+\infty\Big|S_0^{(d)}=x\Big)+b_\lambda,
\]
where $b_\lambda=\frac{\lambda-1-2\lambda H(d)}{\lambda+1}$. According to the Markov property of random walk and direct calculation, it is easy to see that $\{L(x):~x\in \mathbb{Z}^d\}$ is an eigenvector of $G$ with respect to eigenvalue $0$ and hence is an eigenvector of $A_t$ with respect to eigenvalue $e^{0t}=1$, i. e.
\begin{equation}\label{equ 2.10}
\sum_{y\in \mathbb{Z}^d}A_t(x,y)L(y)=L(x)
\end{equation}
for any $x\in \mathbb{Z}^d$. Since $\lim_{d\rightarrow+\infty}H(d)=0$, we assume that $d$ is sufficiently large such that $b_\lambda>0$, then
by Equation \eqref{equ 2.10},
\begin{equation}\label{equ 2.11}
F_t(u)=\sum_{y\in \mathbb{Z}^d}A_t(u,y)\leq \sum_{y\in \mathbb{Z}^d}A_t(u,y)\frac{L(y)}{b_\lambda}=\frac{L(u)}{b_\lambda}.
\end{equation}
For $u\neq 0$, it is obviously that $L(u)\leq H(d)+b_\lambda$, as a result,
\begin{equation}\label{equ 2.12}
\sum\limits_{x\in A}\sum\limits_{y\in A}F_t(x-y)\leq \frac{(|A|^2-|A|)(H(d)+b_\lambda)+|A|(1+b_\lambda)}{b_\lambda}
\end{equation}
according to Equation \eqref{equ 2.11}. By Equations \eqref{equ 2.5}, \eqref{equ 2.9} and \eqref{equ 2.12},
\begin{align}\label{equ 2.13}
&P_{\lambda,d}(\eta_t^{\mathbb{Z}^d}(x)=1 \text{~for some~}x\in A)  \notag\\
&\geq \frac{|A|^2b_\lambda}{(|A|^2-|A|)(H(d)+b_\lambda)+|A|(1+b_\lambda)}\notag\\
&=\frac{|A|^2\big(\lambda-1-2\lambda H(d)\big)}{(|A|^2-|A|)(\lambda-1)\big(1-H(d)\big)+2|A|\lambda\big(1-H(d)\big)}
\end{align}
for any $t>0$. Lemma \ref{lemma 2.1} follows from Equations \eqref{equ 2.4} and \eqref{equ 2.13} directly.

\qed

\quad

\textbf{Acknowledgments.} The author is grateful to the financial
support from the National Natural Science Foundation of China with
grant number 11501542 and the financial support from Beijing
Jiaotong University with grant number KSRC16006536.

{}

\begin{thebibliography}{}
\bibitem{Grif1983}Griffeath, D. (1983). The Binary Contact Path Process. \emph{The Annals of Probability} \textbf{11} 692-705.
\bibitem{Har1974}Harris, T. E. (1974). Contact interactions on a lattice. \emph{The Annals of Probability} \textbf{2}, 969-988.
\bibitem{Hol1981}Holley, R. and Liggett, T. M. (1981). Generalized potlatch and smoothing processes.
\emph{Zeitschrift f\"{u}r Wahrscheinlichkeitstheorie und Verwandte Gebiete } \textbf{55}, 165-195.
\bibitem{Kesten1964}Kesten, H. (1964). On the number of self-avoiding
walks \uppercase\expandafter{\romannumeral2}. \emph{Journal of
Mathematical Physics} \textbf{5}, 1128-1137.
\bibitem{Lig1985}Liggett, T. M. (1985). \emph{Interacting Particle Systems.} Springer, New York.
\bibitem{Lig1999}Liggett, T. M. (1999). \emph{Stochastic interacting systems: contact, voter and exclusion processes.}
Springer, New York.
\end{thebibliography}
\end{document}